\documentclass[12pt]{article}

\usepackage{amsmath, amsthm, amssymb, enumerate}
\usepackage{color}
\usepackage{datetime}
\usepackage{fancyhdr}
\usepackage{eufrak}
\usepackage{enumitem}
\chardef\coloryes=1 
\chardef\isitdraft=1 

\ifnum\isitdraft=1 
   \def\version{9} 
   \def\eqref#1{({\ref{#1}})}                

\definecolor{labelkey}{gray}{.3}

\definecolor{refkey}{rgb}{.3,0.3,0.3}

\else

  \def\startnewsection#1#2{\section{#1}\label{#2}\setcounter{equation}{0}}   
  \def\nnewpage{} 
\fi
\DeclareMathOperator*{\esssup}{ess\,sup}

 \usepackage{amsmath}
\usepackage{amsfonts}
\usepackage{amssymb}
\def\sg {{\sf g}}

\newcommand{\cA}{{\mathfrak A}}

\newcommand{\cP}{{\mathfrak P}}

\newcommand{\cM}{{\mathfrak M}}

\newcommand{\cB}{{\mathfrak B}}

\newtheorem{theorem}{Theorem}[section]
\newtheorem{lemma}{Lemma}[section]

\newtheorem{remark}{Remark}[section]
\newtheorem{definition}{Definition}[section]

\newtheorem{cor}{Corollary}[section]

\newcommand{\bgeqn}{\begin{eqnarray}}
\newcommand{\edeqn}{\end{eqnarray}}
\newcommand{\bgeq}{\begin{eqnarray*}}
\newcommand{\edeq}{\end{eqnarray*}}
\newcommand{\bec}{\begin{center}}
\newcommand{\enc}{\end{center}}
\newcommand{\beqn}{\begin{equation}}
\newcommand{\eeqn}{\end{equation}}
\newcommand{\belign}{\begin{align}}
\newcommand{\elign}{\end{align}}
\newcommand{\belem}{\begin{lemma}}
\newcommand{\elem}{\end{lemma}}
\newcommand{\becor}{\begin{Corollary}}
\newcommand{\ecor}{\end{Corollary}}
\newcommand{\bethm}{\begin{theorem}} 
\newcommand{\ethm}{\end{theorem}}
\newcommand{\beass}{\begin{Assumption}} 
\newcommand{\eass}{\end{Assumption}}
\newcommand{\beitem}{\begin{itemize}}
\newcommand{\eitem}{\end{itemize}}
\newcommand{\beproof}{\begin{proof}}
\newcommand{\eproof}{\end{proof}}
\newcommand{\bedef}{\begin{definition}}
\newcommand{\edefn}{\end{definition}}
\newcommand{\berem}{\begin{remark}}
\newcommand{\erem}{\end{remark}}

\newcommand{\trm}{\textrm}

\newcommand{\Z}{{\cal Z}}
\newcommand{\F}{{\cal F}}

\newcommand{\A}{{\cal A}}
\newcommand{\B}{{\cal B}}

\newcommand{\M}{{\cal M}}

\newcommand{\bbP}{\mathbb{P}}
\newcommand{\bbN}{\mathbb{N}_0}
\newcommand{\bbR}{\mathbb{R}}
\newcommand{\bbK}{\mathbb{K}}
\newcommand{\bbF}{\mathbb{F}}

\newcommand{\ra}{\rightarrow}

\newcommand{\be}{\begin{equation}}
\newcommand{\ee}{\end{equation}}
\newcommand{\rar}{\rightarrow}

\newcommand{\uar}{\uparrow}

\newcommand{\bbpara}{\bigg\{ }
\newcommand{\bbpark}{\bigg\} }
\newcommand{\la}{\langle }
\newcommand{\ran}{\rangle }

\newcommand{\varo}{{\varrho}}

\newcommand{\lan}{\langle}

\def\beal#1\eal{\begin{align}#1\end{align}}
\def\dist{\mathop{\rm dist}}
\def\w{\omega}

\def\O{\Omega}

\def\ra{{\rangle}}

\def\supp {{\rm supp}}

\newcommand{\avr}{{\sf AV@R}}

\newcommand{\VaR}{{\sf V@R}}
\newcommand{\essinf}{\operatornamewithlimits{ess\,inf}}

\def\bbp{{\Bbb{P}}}

\def\bbP{{\mathbb{P}}}

\def\bbE{{\mathbb{E}}}
\def \rhot1{\rho_{t+1}}
\def \rhot{\rho_{t}}
\def \alf{\alpha}
\def \sg{\sigma}

 \renewcommand{\theequation}{\thesection.\arabic{equation}}

\begin{document}
\def\ques{{\cor \underline{??????}\cob}}
\def\nto#1{{\coC \footnote{\em \coC #1}}}
\def\fractext#1#2{{#1}/{#2}}
\def\fracsm#1#2{{\textstyle{\frac{#1}{#2}}}}   
\def\nnonumber{}


\def\cor{{}}
\def\cog{{}}
\def\cob{{}}
\def\coe{{}}
\def\coA{{}}
\def\coB{{}}
\def\coC{{}}
\def\coD{{}}
\def\coE{{}}
\def\coF{{}}
\ifnum\coloryes=1

  \definecolor{coloraaaa}{rgb}{0.1,0.2,0.8}
  \definecolor{colorbbbb}{rgb}{0.1,0.7,0.1}
  \definecolor{colorcccc}{rgb}{0.8,0.3,0.9}
  \definecolor{colordddd}{rgb}{0.0,.5,0.0}
  \definecolor{coloreeee}{rgb}{0.8,0.3,0.9}
  \definecolor{colorffff}{rgb}{0.8,0.3,0.9}
  \definecolor{colorgggg}{rgb}{0.5,0.0,0.4}

 \def\cog{\color{colordddd}}
 \def\cob{\color{black}}
 \def\cor{\color{red}}
 \def\coe{\color{colorgggg}}

 \def\coA{\color{coloraaaa}}
 \def\coB{\color{colorbbbb}}
 \def\coC{\color{colorcccc}}
 \def\coD{\color{colordddd}}
 \def\coE{\color{coloreeee}}
 \def\coF{\color{colorffff}}
 \def\coG{\color{colorgggg}}

\fi
\ifnum\isitdraft=1
   \chardef\coloryes=1 
   \baselineskip=17pt
\pagestyle{myheadings}
\reversemarginpar

\def\const{\mathop{\rm const}\nolimits}  
\def\diam{\mathop{\rm diam}\nolimits}    

 \def\llabel#1{\label{#1}{\ \mbox{\rm\color{red} {#1}\color{black}}}}

\def\rref#1{{\ref{#1}{\rm \tiny \fbox{\tiny #1}}}}
\def\theequation{\fbox{\bf \thesection.\arabic{equation}}}
\def\ccite#1{{\cite{#1}{\rm \tiny ({#1})}}}

\def\startnewsection#1#2{\newpage\cog \section{#1}\cob\label{#2}

\setcounter{equation}{0}
\pagestyle{fancy}

\lhead{\cob Section~\ref{#2}, #1 }
\cfoot{}
\rfoot{\thepage}
\lfoot{\cob{\today,~\currenttime}~(c75-iklt2, Version~\fbox{\version})}}
\chead{}
\rhead{\thepage}
\def\nnewpage{\newpage}

\newcounter{startcurrpage}
\newcounter{currpage}

\def\llll#1{{\rm\tiny\fbox{#1}}}
   \def\blackdot{{\color{red}{\hskip-.0truecm\rule[-1mm]{4mm}{4mm}\hskip.2truecm}}\hskip-.3truecm}
   \def\bdot{{\coC {\hskip-.0truecm\rule[-1mm]{4mm}{4mm}\hskip.2truecm}}\hskip-.3truecm}
   \def\purpledot{{\coA{\rule[0mm]{4mm}{4mm}}\cob}}
   \def\pdot{\purpledot}
\else
   \baselineskip=15pt
   \def\blackdot{{\rule[-3mm]{8mm}{8mm}}}
   \def\purpledot{{\rule[-3mm]{8mm}{8mm}}}
   \def\pdot{}
\fi

\def\tdot{\fbox{\fbox{\bf\tiny I'm here; \today \ \currenttime}}}
\def\nts#1{{\hbox{\bf ~#1~}}} 
\def\nts#1{{\cor\hbox{\bf ~#1~}}} 
\def\ntsf#1{\footnote{\hbox{\bf ~#1~}}} 
\def\ntsf#1{\footnote{\cor\hbox{\bf ~#1~}}} 
\def\bigline#1{~\\\hskip2truecm~~~~{#1}{#1}{#1}{#1}{#1}{#1}{#1}{#1}{#1}{#1}{#1}{#1}{#1}{#1}{#1}{#1}{#1}{#1}{#1}{#1}{#1}\\}
\def\biglineb{\bigline{$\downarrow\,$ $\downarrow\,$}}
\def\biglinem{\bigline{---}}
\def\biglinee{\bigline{$\uparrow\,$ $\uparrow\,$}}

\def\w{\omega}
\def\tilde{\widetilde}

\newtheorem{Theorem}{Theorem}[section]
\newtheorem{Corollary}[Theorem]{Corollary}
\newtheorem{Proposition}[Theorem]{Proposition}
\newtheorem{Lemma}[Theorem]{Lemma}
\newtheorem{Remark}[Theorem]{Remark}
\newtheorem{Example}[Theorem]{Example}
\newtheorem{Assumption}[Theorem]{Assumption}
\newtheorem{Claim}[Theorem]{Claim}
\newtheorem{Question}[Theorem]{Question}
\def\theequation{\thesection.\arabic{equation}}
\def\endproof{\hfill$\Box$\\}
\def\square{\hfill$\Box$\\}
\def\comma{ {\rm ,\qquad{}} }            
\def\commaone{ {\rm ,\qquad{}} }         
\def\dist{\mathop{\rm dist}\nolimits}    
\def\sgn{\mathop{\rm sgn\,}\nolimits}    
\def\Tr{\mathop{\rm Tr}\nolimits}    
\def\div{\mathop{\rm div}\nolimits}    
\def\supp{\mathop{\rm supp}\nolimits}    
\def\divtwo{\mathop{{\rm div}_2\,}\nolimits}    
\def\re{\mathop{\rm {\mathbb R}e}\nolimits}    
\def\div{\mathop{\rm{Lip}}\nolimits}   
\def\indeq{\qquad{}}                     
\def\period{.}                           
\def\semicolon{\,;} 


\title{Robust Optimal Control Using Conditional Risk Mappings in Infinite Horizon}
\author{Kerem U\u{g}urlu}
\maketitle

\date{}

\begin{center}
\end{center}

\medskip

\indent Department of Applied Mathematics, University of Washington, Seattle, WA 98195\\
\indent e-mail: keremu@uw.edu

\begin{abstract}
We use one-step conditional risk mappings to formulate a risk averse version of a total cost problem on a controlled Markov process in discrete time infinite horizon. The nonnegative one step costs are assumed to be lower semi-continuous but not necessarily bounded. We derive the conditions for the existence of the optimal strategies and solve the problem explicitly by giving the robust dynamic programming equations under very mild conditions. We further give an $\epsilon$-optimal approximation to the solution and illustrate our algorithm in two examples of optimal investment and LQ regulator problems. 
\end{abstract}

\section{Introduction}
Controlled Markov decision processes have been an active research area in sequential decision making problems in operations research and in mathematical finance. We refer the reader to \cite{HL,HL2,SB} for an extensive treatment on theoretical background. Classically, the evaluation operator has been the expectation operator, and the optimal control problem is to be solved via Bellman's dynamic programming \cite{key-6}. This approach and the corresponding problems continue to be an active research area in various scenarios (see e.g. the recent works \cite{SN,NST,GHL} and the references therein) 

On the other hand, expected values are not appropriate to measure the performance of the agent. Hence, expected criteria with utility functions have been extensively used in the literature (see e.g. \cite{FS1,FS2} and the references therein). Other than the evaluation of the performance via utility functions, to put risk aversion into an axiomatic framework, coherent risk measures has been introduced in the seminal paper \cite{key-1}. \cite{FSH1} has removed the positive homogeneity assumption of a coherent risk measure and named it as a convex risk measure (see \cite{FSH2} for an extensive treatment on this subject). 

However, this kind of operator has brought up another difficulty.  Deriving dynamic programming equations with these operators in multistage optimization problems is challenging or impossible in many optimization problems. The reason for it is that the Bellman's optimality principle is not necessarily true using this type of operators. That is to say, the optimization problems are not \textit{time-consistent}. Namely, a multistage stochastic decision problem is time-consistent, if resolving the problem at later stages (i.e., after observing some random outcomes), the original solutions remain optimal for the later stages. We refer the reader to \cite{key-141, key-334, key-148, AS, BMM} for further elaboration and examples on this type of inconsistency. Hence, optimal control problems on multi-period setting using risk measures on bounded and unbounded costs are not vast, but still, some works in this direction are \cite{ key-149,key-150,key-151, BO}. 

To overcome this deficit, dynamic extensions of convex/coherent risk measures so called conditional risk measures are introduced in \cite{FR} and studied extensively in \cite{RS06}. In \cite{key-3}, so called Markov risk measures are introduced and an optimization problem is solved in a controlled Markov decision framework both in finite and discounted infinite horizon, where the cost functions are assumed to be bounded. This idea is extended to transient models in \cite{CR1,CR2} and to unbounded costs with $w$-weighted bounds in \cite{LS,CZ,SSO} and to so called \textit{process-based} measures in \cite{FR1} and to partially observable Markov chain frameworks in \cite{FR2}.

In this paper, we derive \textit{robust} dynamic programming equations in discrete time on infinite horizon using one step conditional risk mappings that are dynamic analogues of coherent risk measures. We assume that our one step costs are nonnegative, but may well be unbounded from above. We show the existence of an optimal policy via dynamic programming under very mild assumptions. Since our methodology is based on dynamic programming, our optimal policy is by construction time consistent. We further give a recipe to construct an $\epsilon$-optimal policy for the infinite horizon problem and illustrate our theory in two examples of optimal investment and LQ regulator control problem, respectively. To the best of our knowledge, this is the first work solving the optimal control problem in infinite horizon with the minimal assumptions stated in our model.

The rest of the paper is as follows. In Section 2, we briefly review the theoretical background on coherent risk measures and their dynamic analogues in multistage setting, and further describe the framework for the controlled Markov chain that we will work on. In Section 3, we state our main result on the existence of the optimal policy and the existence of optimality equations. In Section 4, we prove our main theorem and present an $\epsilon$ algorithm to our control problem. In Section 5, we illustrate our results with two examples, one on an optimal investment problem, and the other on an LQ regulator control problem.
\section{Theoretical Background}

In this section, we recall the necessary background on static coherent risk measures, and then we extend this kind of operators to the dynamic setting in controlled Markov chain framework in discrete time. 

\subsection{Coherent Risk Measures}
Consider an atomless probability space $(\O,\F, \bbP)$ and the space $\Z: = L^1(\O,\F,\bbP)$ of measurable functions $Z:\O \rar \bbR$ (random variables) having finite first order moment, i.e. $\bbE^{\bbP}[|Z|] < \infty$, where $\bbE^\bbP[\cdot]$ stands for the expectation with respect to the probability measure $\bbP$. A mapping $\rho:\Z \rar \bbR$ is said to
be a \textit{coherent risk measure}, if it satisfies the following axioms
\begin{itemize}
\item (A1)(Convexity) $\rho(\lambda X+(1-\lambda)Y)\leq\lambda\rho(X)+(1-\lambda)\rho(Y)$
$\forall\lambda\in(0,1)$, $X,Y \in \Z$.
\item (A2)(Monotonicity) If $X \preceq Y$, then $\rho(X) \leq \rho(Y)$, for all $X,Y \in \Z$.
\item (A3)(Translation Invariance) $\rho(c+X) = c + \rho(X)$, $\forall c\in\mathbb{R}$, $X\in \Z$.
\item (A4)(Homogeneity) $\rho(\beta X)=\beta\rho(X),$ $\forall X\in \Z$. $\beta \geq 0$.
\end{itemize}
The notation $X \preceq Y$ means that $X(\w) \leq Y(\w)$ for $\bbP$-a.s. Risk measures $\rho:\Z \rar \bbR$, which satisfy (A1)-(A3) only, are called convex risk measures. We remark that under the fourth property (homogeneity), the first property (convexity) is equivalent to sub-additivity. We call the risk measure $\rho:\Z\rar \bbR$ law invariant, if $\rho(X) = \rho(Y)$, whenever $X$ and $Y$ have the same distributions.  We pair the space $\Z = L^1(\O,\F, \bbP)$ with $\Z^* = L^\infty(\O,\F, \bbP)$, and the corresponding scalar product 
\beqn 
\lan \zeta, Z \ran = \int_\O \zeta(\w)Z(\w)dP(\w), \; \zeta\in \Z^*,Z\in\Z.
\eeqn
By \cite{key-14}, we know that real-valued law-invariant convex risk measures are continuous, hence lower semi-continuous (l.s.c.), in the norm topology of the space $L^1(\O,\F, \bbP)$. 
Hence, it follows by Fenchel-Moreau theorem that
\beqn 
\label{eqn12}
\rho(Z) = \sup_{\zeta \in \Z^*} \{ \lan \zeta, Z \ran - \rho^*(\zeta) \},\;\textrm{for all } Z \in \Z,
\eeqn 
where $\rho^*(Z) = \sup_{Z \in \Z} \{ \lan \zeta, Z \ran - \rho(Z) \}$ is the corresponding conjugate functional (see \cite{RW}). If the risk measure $\rho$ is convex and positively homogeneous, hence coherent, then $\rho^*$ is an indicator function of a convex and closed set $\cA \subset \Z^*$ in the respective paired topology. The dual representation in Equation \ref{eqn12} then takes the form
\beqn 
\label{eqn13}
\rho(Z) = \sup_{\zeta \in \cA} \lan \zeta,Z \ran,\;Z\in \Z,
\eeqn 
where the set $\cA$ consists of probability density functions $\zeta:\O\rar \bbR$, i.e. with $\zeta \succeq 0$ and $\int\zeta dP = 1$. 

A fundamental example of law invariant coherent risk measures is Average- Value-at-Risk measure (also
called the Conditional-Value-at-Risk or Expected Shortfall Measure). Average-Value- at-Risk at the level of $\alpha$ for $Z \in \Z$ is defined as
\beqn
\label{eqn016}
\avr_\alpha (Z) = \frac{1}{1-\alpha}\int_\alpha^1 \VaR_p(Z)dp,
\eeqn 
where
\beqn 
\VaR_p(Z) = \inf \{ z \in \bbR: \bbP(Z \leq z) \geq p \}
\eeqn 
is the corresponding left side quantile. The corresponding dual representation for $\avr_\alpha(Z)$ is 
\beqn
\avr_\alpha(Z) = \sup_{m \in \A}\lan m,Z \ran,
\eeqn
with
\beqn 
\A = \{ m \in L^\infty(\O,\F,\bbP): \int_\O md\bbP =1, 0 \leq \lVert m \rVert_\infty \leq \frac{1}{\alpha} \}.
\eeqn 
Next, we give a representation characterizing any law invariant coherent risk measure, which is first presented in Kusuoka \cite{K} for random variables in  $L^\infty(\O, \F, \bbP)$, and later further investigated in $\Z^p = L^p(\O, \F, \bbP)$ for $1 \leq p < \infty$ in \cite{PR}.

\belem
\label{lem11}\cite{K}
Any law invariant coherent risk measure $\rho: \Z^p \rar \bbR$ can be represented in the following
form
\beqn 
\rho(Z) = \sup_{\nu \in \cM}\int_0^1 \avr_\alpha(Z)d\nu(\alpha),
\eeqn 
where $\cM$ is a set of probability measures on the interval [0,1].
\elem 

\subsection{Controlled Markov Chain Framework}
Next, we introduce the controlled Markov chain framework that we are going to study our problem on. 
We take the control model $\mathcal{M} = \{ \mathcal{M}_n, n \in \mathbb{N}_0 \}$, where for each $n \geq 0$, we have 

\beqn 
\label{eqn270}
\M_n := (X_n, A_n, \mathbb{K}_n, Q_n, F_n, c_n)
\eeqn 

with the following components:

\beitem 
\item $X_n$ and $A_n$ denote the state and action (or control) spaces,which are assumed to be complete seperable metric spaces with their corresponding Borel $\sg$-algebras $\B(X_n)$ and $\B(A_n)$. 
\item For each $x_n \in X_n$, let $A_n(x_n) \subset A_n$ be the set of all admissible controls in the state $x_n$.  Then 
\begin{equation}
\mathbb{K}_n := \{ (x_n,a_n): x_n \in X_n,\; a_n \in A_n \}
\end{equation}
stands for the set of feasible state-action pairs at time $n$.

\item We let 
\begin{equation}
\label{eqn23}
x_{i+1} = F_i(x_i,a_i, \xi_i), 
\end{equation}
for all $i = 0,1,...$ with $x_i \in X_i$ and $a_i \in A_i$ as described above, with independent random variables $(\xi_i)_{i \geq 0}$ on the atomless probability space 
\beqn 
\label{eqn2100}
(\O^i,\mathcal{G}^i, \bbP^i). 
\eeqn 
We take that $\xi_i \in S_i$, where $S_i$ are Borel spaces.  Moreover, we assume that the system equation  
\beqn 
\label{eqn678}
F_i: \mathbb{K}_i \times S_i \rar X_i
\eeqn 
as in Equation \eqref{eqn23} is continuous. 

\item We let 
\beal
\label{eqn27}
\O &= \otimes_{i=1}^{\infty} X^i
\eal 
where $X^i$ is as defined in Equation \eqref{eqn678}. For $n\geq 0$, we let 
\beal
\label{eqn28}
\F_n &= \sg(\sg({\displaystyle \cup_{i=0}^n  \mathcal{G}^i)} \cup \sg(X_0,A_0,X_1,A_1\ldots,A_{n-1},X_n)) \\
\F &= \sg(\cup_{i=0}^\infty \F_{i})
\eal be the filtration of increasing $\sg$-algebras. Furthermore, we define the corresponding probability measures $(\O,\F)$ as 
\beal
\label{eqn2140}
\bbP &= \prod_{i=1}^{\infty} \bbp^i,
\eal
where the existence of $\bbP$ is justified by Kolmogorov extension theorem (see \cite{HL}). We assume that for any $n\geq 0$, the random vector $\xi_{[n]} = (\xi_0,\xi_1,\ldots,\xi_n)$ and $\xi_{n+1}$ are independent on $(\O,\F,\bbP)$. 

\item The transition law is denoted by $Q_{n+1}(B_{n+1}|x_n,a_n)$, where $B_{n+1} \in \mathcal{B}(X_{n+1})$ is the Borel $\sigma$-algebra on $X_n$, and $(x_n,a_n) \in X_n \times A_n$ is a stochastic kernel on $X_n$ given $\mathbb{K}_n$ (see \cite{SB,HL} for further details). We remark here that at each $n\geq 0$ the stochastic kernel depends only on $(x_n,a_n)$ rather than $\F_n$.
That is, for each pair $(x_n,a_n) \in \mathbb{K}_n$, $Q_{n+1}(\cdot|x_n,a_n)$ is  a probability measure on $X_{n+1}$, and for each $B_{n+1} \in \mathcal{B}_{n+1}(X_{n+1})$, $Q_{n+1}(B_{n+1}|\cdot,\cdot)$ is a measurable function on $\mathbb{K}_n$. Let  $x_0 \in X_0$ be given with the corresponding policy $\Pi = (\pi_n)_{n\geq0}$. By the Ionescu Tulcea theorem (see e.g. \cite{HL}), we know that there exists a unique probability measure $\bbP^\pi$ on $(\O, \F)$ such that given $x_0 \in X_0$, a measurable
set $B_{n+1} \subset X_{n+1}$ and $(x_n,a_n) \in \mathbb{K}_n$, for any $n\geq 0$, we have 
\beal
\label{eqn2150}
\bbP^{\Pi}_{n+1} (x_{n+1} \in B_{n+1}) &\triangleq Q_{n+1}(B_{n+1} | x_n,a_n).
\eal
\item Let $\bbF_n$ be the family of measurable functions  $\pi_n:X_n \rar A_n$ for $n \geq 0$. A sequence $( \pi_n )_{n\geq 0}$ of functions $\pi_n \in \bbF_n$ for $n \geq 0$ is called a control policy (or simply a policy), and the function $\pi_n(\cdot)$ is called the decision rule or control at time $n\geq 0$. We denote by $\Pi$ the set of all control policies. For notational convenience, for every $n \in \bbN$ and $(\pi_n)_{n \geq 0} \in \Pi$, we write 
\begin{align*}
c_n(x_n,\pi_n) &:= c_n(x_n,\pi_n(x_n))\\
&:= c_n(x_n,a_n).
\end{align*}
We denote by $\cP(A_n(x_n))$ as the set of probability measures on $A_n(x_n)$ for each time $n\geq0$. A randomized Markovian policy $(\pi_n)_{n \geq 0}$ is a sequence of measurable functions such that $\pi_n(x_n) \in \cP(A_n(x_n))$ for all $x_n \in X_n$, i.e. $\pi_n(x_n)$ is a probability measure on $A_n(x_n)$. $(\pi_n)_{n \geq 0}$ is called a deterministic policy, if $\pi_n(x_n) = a_n$ with $a_n \in A_n(x_n)$.

\item $c_n(x_n,a_n): \mathbb{K}_n \rightarrow \mathbb{R}_{+}$ is the real-valued cost-per-stage function at stage $n \in \mathbb{N}_0$ with  $(x_n,a_n) \in \mathbb{K}_n$.

\eitem

\bedef
\label{defn31}
A real valued function $v$ on $\mathbb{K}_n$ is said to be inf-compact on $\mathbb{K}_n$, if the set 
\begin{equation}
\{ a_n \in A_n(x_n) | v(x_n,a_n) \leq r \}
\end{equation}
is compact for every $x_n \in X_n$ and $r \in \mathbb{R}$. As an example, if the sets $A_n(x_n)$ are compact and $v(x_n,a_n)$ is l.s.c. in $a_n\in A_n(x_n)$ for every $x_n \in X_n$, then $v(\cdot,\cdot)$ is inf-compact on $\mathbb{K}_n$. Conversely, if $v$ is inf-compact on $\mathbb{K}_n$, then $v$ is l.s.c. in $a_n \in A_n(x_n)$ for every $x_n \in X_n$. 
\edefn
We make the following assumption about the transition law $(Q_n)_{n\geq1}$.
\begin{Assumption}\label{ass31}
For any $n \geq 0$, the transition law $Q_n$ is weakly continuous; i.e. for any continuous and bounded function $u(\cdot)$ on $X_{n+1}$, the map 
\begin{equation}
(x_n,a_n) \rightarrow \int_{X_{n+1}} u(y)dQ_n(y|x_n,a_n) 
\end{equation}
is continuous on $\mathbb{K}_n$.
\end{Assumption}

Furthermore, we make the following assumptions on the one step cost functions and action sets.

\beass \label{ass32}
For every $n \geq 0$,
\beitem
\item the real valued non-negative cost function $c_n(\cdot,\cdot)$ is l.s.c. in $(x_n,a_n)$. That is for any $(x_n,a_n) \in X_n \times A_n$, we have
\beqn 
c_n(x_n,a_n) \leq \liminf_{(x^k_n,a^k_n) \rar (x_n,a_n) } c_n(x^k_n,a^k_n),
\eeqn 
as $k\rar \infty$.

\item The multifunction (also known as a correspondence or point-to-set function) $x_n \rar A_n(x_n)$, from $X_n$ to $A_n$, is upper semicontinuous (u.s.c.) that is, if $\{x_n^l\} \subset X_n$ and $\{ a_n^l\} \subset A_n$ are sequences such that $\{x_n^l\} \rar \bar{x}_n$ with $\{a_n^l\} \subset A_n$ for all $l$, and $a_n^l \rar \bar{a}_n$, then $\bar{a}_n$ is in $A_n(\bar{x}_n)$. 

\item For every state $x_n \in X_n$, the admissible action set $A_n(x_n)$ is compact.
\eitem 
\eass

\subsection{Conditional Risk Mappings} In order to construct dynamic models of risk, we extend the concept of static coherent risk measures to dynamic setting. For any $n \geq 1$, we denote the space $\Z_n: = L^1(\O,\F_n,\bbP^\pi_{n})$ of measurable functions with  $Z:\O \rar \bbR$ (random variables) having finite first order moment, i.e. $\bbE^{\bbP^\pi_{n}}[|Z|] < \infty$ $\bbP^\pi_{n}$-a.s., where $\bbE^{\bbP^\pi_{n}}$ stands for the conditional expectation at time $n$ with respect to the conditional probability measure $\bbP^\pi_{n}$ as defined in Equation \eqref{eqn2150}. 
\bedef\label{def21}
Let $X,Y \in \Z_{n+1}$. We say that a mapping $\rho_n:\Z_{n+1} \rar \Z_n$ is a one step conditional risk mapping, if it satisfies following properties
\beitem
\item(a1) Let $\gamma \in [0,1]$. Then, 
\beqn 
\rho_{n}(\gamma X + (1-\gamma) Y) \preceq
\gamma \rho_{n}(X) + (1-\gamma)\rho_{n}(Y)
\eeqn 
\item (a2) If $X \preceq Y$, then $\rho_{n}(X) \preceq \rho_{n}(Y)$
\item (a3) If $Y \in \Z_{n}$ and $X \in \Z_{n+1}$, then $\rho_{n}(X + Y) = \rho_{n}(X) + Y$.
\item (a4) For $\lambda \succeq 0$ with $\lambda \in \Z_n$ and $X \in \Z_{n+1}$, we have that $\rho_{n+1}(\lambda X) = \lambda \rho_{n+1}(X)$.  
\eitem
\edefn
Here, the relation $Y(\w) \preceq X(\w)$ stands for $Y \leq X$ $\bbP^{\pi}_n$-a.s. We next state the analogous results for representation theorem for conditional risk mappings as in Equation \eqref{eqn13} (see also \cite{RS06}).

\bethm
\label{thm21}
Let $\rho_n: \Z_{n+1} \rar \Z_n$ be a law-invariant conditional risk mapping satisfying assumptions as stated in Definition \ref{def21}. Let $Z \in \Z_{n+1}$. Then
\beqn 
\label{eqn2220}
\rho_n(Z) = \sup_{\mu \in \cA_{n+1}}  \lan \mu, Z \ran,
\eeqn 
where $\cA_{n+1}$ is a convex closed set of conditional probability measures on  $(\O, \F_{n+1})$, that are absolutely continuous with respect to  $\bbP^{\pi}_{n+1}$. 
\ethm
Next, we give the Kusuoka representation for conditional risk mappings analogous to Lemma \ref{lem11}. 
\belem \label{lem22}
Let $\rho_n:\Z_{n+1} \rar \Z_{n}$ be a law invariant one-step conditional risk mapping satisfying Assumptions
(a1)-(a4) as in Definition \ref{def21}. Let $Z \in \Z_{n+1}$. Then, conditional Average-Value-at-Risk at the level of $0 < \alpha < 1$ is defined as 
\beqn
\label{eqn016}
\avr^{n}_\alpha (Z) \triangleq \frac{1}{1-\alpha}\int_\alpha^1 \VaR^{n}_p(Z)dp,
\eeqn 
where 
\beqn 
\label{eqn212}
\VaR^{n}_p(Z) \triangleq \essinf \{ z \in \bbR: \bbP^{\pi}_{n+1}(Z \leq z) \geq p \}.
\eeqn
Here, we note that $\VaR^{n}_p(Z)$ is $\F_{n}$-measurable by definition of essential infimum (see \cite{FSH2} for a definition of essential infimum and essential supremum). 
Then, we have 
\beqn 
\label{eqn2130}
\rho_{n}(Z) \triangleq  \esssup_{\nu \in \cM}\int_0^1\avr^n_\alpha(Z)d\nu(\alpha),
\eeqn 
where $\cM$ is a set of probability measures on the interval [0,1].  
\elem 

\berem By Equations \eqref{eqn016},\eqref{eqn212} and \eqref{eqn2130}, it is easy to see that the corresponding optimal controls at each time $n \geq 0$ is deterministic, if the one step conditional risk mappings are $\avr^n_\alf: \Z_{n+1} \rar \Z_{n}$ as defined in \eqref{eqn016}. On the other hand, by Kusuoka representation, Equation \eqref{eqn2130}, it is clear that for other  coherent risk randomized policies might be optimal. In this paper, we restrict our study to deterministic policies.
\erem

\bedef\label{defn230} A policy $\pi \in \Pi$ is called admissible, if for any $n\geq 0$, we have 
\beal 
&c_n(x_n,a_n) + \lim_{N\rar \infty}\gamma\rho_n\big( c_{n+1}(x_{n+1},a_{n+1}) \\
&\indeq + \gamma\rho_{n+1}( c_{n+2}(x_{n+2},a_{n+2}) 
\ldots  + \gamma\rho_{N-1}(c_N(x_N,a_N))) \big) < \infty,\; \bbP^\pi_{n}\trm{ a.s. }
\eal
The set of all admissible policies is denoted by $\Pi_{\mathrm{ad}}$. 
\edefn

\section{Main Problem}

Under Assumptions \ref{ass31}, \ref{ass32}, our control problem reads as 
\beal 
\label{eqn321}
&\inf_{\pi \in \Pi_{\textrm{ad}}} \bigg( c_0(x_0,a_0) + \lim_{N\rar \infty}\gamma\rho_0( c_1(x_1,a_1) + \gamma\rho_{1}( c_2(x_2,a_2) \\
&\indeq \ldots  + \gamma\rho_{N-1}(c_N(x_N,a_N))) \bigg)
\eal
Namely, our objective is to find a policy $(\pi^*_n)_{ n \geq 0}$ such that the value function in Equation \eqref{eqn321} is minimized. For convenience, we introduce the following notations that are to be used in the rest of the paper
\begin{align*}
\varrho_{n-1} ( \sum_{t=n}^\infty c_t(x_t,\pi_t) ) &:= \lim_{N\rar \infty}\gamma\rho_{n-1} (c_n(x_n,a_n) + \gamma\rho_{n}( c_{n+1}(x_{n+1},a_{n+1}) 
\\
&...  + \rho_{N-1}(c_N(x_N,a_N))) \\
V_n(x,\pi) &:= c_n(x_n,a_n) + \varrho_n(\sum_{t=n+1}^\infty c_t(x_t,a_t))
\end{align*}
\begin{align*}
V_n^*(x) &:= \inf_{\pi \in \Pi_{\trm{ad}}} c_n(x_n,a_n) + \varrho_n(\sum_{t=n+1}^\infty c_t(x_t,a_t))\\
V_{n,N}(x,\pi) &:= c_n(x_n,a_n) + \varrho_n(\sum_{t=n+1}^{N-1}c_t(x_t,a_t)) \\
V_{N,\infty}(x,\pi) &:= c_N(x_N,a_N) + \varrho_N(\sum_{t=N+1}^\infty c_t(x_t,a_t)) \\  
V_{n,N}^*(x) &:=  \inf_{\pi \in \Pi_{\trm{ad}}} c_N(x_N,a_N) +\varrho_n(\sum_{t=n+1}^Nc_t(x_t,a_t))
\end{align*}
For the control problem to be nontrivial, we need the following assumption on the existence of the policy.
\begin{Assumption}
\label{ass41}
There exists a policy $\pi \in \Pi_{\trm{ad}}$ such that 
\beqn 
c_0(x_0, a_0) + \varrho_0(x_0)  < \infty.
\eeqn 
\end{Assumption}
We are now ready to state our main theorem. 
\bethm
\label{thm41}
Let $0 < \gamma < 1$. Suppose that Assumptions \ref{ass31}, \ref{ass32} and \ref{ass41} are satisfied. Then,
\beitem
\item[(a)] the optimal cost functions $V_n^*$ are the pointwise minimal solutions of the optimality equations: that is, for every $n \in \bbN$ and $x_n \in X_n$, 
\beqn
\label{eqn33}
V_n^*(x_n) = \inf_{a \in A(x_n)}  \bigg(c_n(x_n,a_n) +  \gamma\rho_n(V_{n+1}^*(x_{n+1}))\bigg).
\eeqn 
\item[(b)] There exists a policy $\pi^* = (\pi^*_n)_{n \geq 0}$ such that for each $n \geq 0$, the control attains the minimum in \eqref{eqn33}, namely for $x_n \in X_n$
\beqn 
V_n^*(x_n) = c_n(x_n,\pi^*_n) +  \gamma\rho_n(V_{n+1}^*(x_{n+1}) ).
\eeqn 
\eitem
\ethm
\section{Proof of Main Result}

\belem \cite{key-35} \label{lem31} Fix an arbitrary $n \in \bbN$. Let $\bbK$ be defined as 
\beqn
\bbK :=  \{ (x,a)| x \in X, a \in A(x) \},
\eeqn 
where $X$ and $A$ are complete seperable metric Borel spaces and let $v: \bbK \rar \bbR$ be a given $\B(X \times A)$ measurable function. For $x \in X$, define 
\beqn 
v^*(x) := \inf_{a \in A(x)} v(x,a).
\eeqn
If $v$ is non-negative, l.s.c. and inf-compact on $\bbK$ as defined in Definition \ref{defn31}, then for any $x \in X$, there exists a measurable mapping $\pi_n: X \rar A$ such that 
\beqn 
v^*(x) = v(x,\pi_n)
\eeqn
and $v^*(\cdot):X\rar \bbR$ is measurable, and l.s.c. 
\elem 


\belem \label{lem52}
For any $n \geq 1$, let $c_n(x_n,a_n)$ be in $\Z_n$. Then $\rho_{n-1}(c_n(x_n,a_n))$ is an element of $\Z_{n-1} = L^1(\O,\F_n,\bbP^\pi_n)$.
\elem 
\beproof
Let $\mu \in \cA_{n}$ be as in Theorem \ref{thm21}. By non-negativity of the one step cost function $c_n(\cdot, \cdot)$ and by Fatou Lemma, we have
\beqn 
\label{eqn433}
\lan \mu,c_n(x_n,a_n) \ra \leq \liminf_{(x^k_n,a^k_n) \rar (x_n,a_n)} \lan \mu,c_n(x^k_n,a^k_n) \ra. 
\eeqn
Hence, $\lan \mu,c_n(x_n,a_n) \ra$ is l.s.c. for  $\bbP^{\pi}_{n-1}$-a.s. Then, by Equation \eqref{eqn2220}, we have 
\beqn 
\label{eqn434}
\rho_{n-1}(c_n(x_n,a_n)) = \esssup_{\mu \in \cA_{n}} \lan \mu, c_n(x_n,a_n) \ran. 
\eeqn 
Hence, by Equation \eqref{eqn433} and by Equation \eqref{eqn434} taking supremum of l.s.c. functions being still l.s.c., we conclude that for fixed $\w$, $\rho_{n-1}(c_n(x_n(\w),a_n(\w)))$ is l.s.c. with respect to $(x_n,a_n)$. 

Next, we show that $\rho_{n-1}(c_n(x_n,a_n))$ is $\F_{n-1}$ measurable. By Lemma \ref{lem22}, we have 
\beal 
\label{eqn018}
&\rho_{n-1}(c_n(x_{n},a_{n})) = \esssup_{\nu \in \cM} \int_{[0,1]} \avr^{n-1}_\alpha(c_n(x_{n},a_{n})) d\nu,\\
&\indeq= \esssup_{\nu \in \cM} \int_{[0,1]} \frac{1}{1-\alpha}\int_\alpha^1 \VaR_p^{n-1}(c_n(x_{n},a_{n})) dp\; d\nu\\
&\indeq= \esssup_{\nu \in \cM} \int_{[0,1]} \frac{1}{1-\alpha}\int_\alpha^1 \essinf \big( z \in \bbR: \bbP^\pi_n(c_n(x_{n},a_{n})\leq z) \geq p \big) dp\; d\nu,
\eal
where $\cM$ is a set of probability measures on the interval [0,1]. By noting that for any $p \in [\alpha,1]$, $\essinf \big( z \in \bbR: \bbP^\pi_n(c_n(x_{n},a_{n})\leq z) \geq p \big)$ is $\F_{n-1}$-measurable, and then, by integrating from $\alpha$ to 1 and multiplying by $\frac{1}{1-\alpha}$, $\F_{n-1}$ measurability is preserved. Similarly, in Equation \ref{eqn018}, integrating with respect to a probability measure $\nu$ on $[0,1]$ and taking supremum of the integrals preserve $\F_{n-1}$ measurability. Hence, we conclude the proof. 
\eproof

\becor Let $n \geq 1$, $x_n \in X_n$ and $a_n \in A_n$, where $X_n$ and $A_n$ are as introduced in Equation \eqref{eqn270}. Then, 
\label{cor51}
\beqn 
\min_{a_n \in \pi(x_n)} \rho_{n-1}(c_n(x_n,a_n))
\eeqn 
is  l.s.c. in $x_n$ $\bbP^{\pi}_{n-1}$-a.s. Furthermore, $\displaystyle\min_{a_n \in \pi(x_n)} \rho_{n-1}(c_n(x_n,a_n))$ is $\F_{n-1}$ measurable. 
\ecor 

\beproof
We know by Lemma \ref{lem52}, $\rho_{n-1}(c_n(x_n,a_n))$ is l.s.c.  $\bbP^{\pi}_{n-1}$-a.s.  Hence, by Lemma \ref{lem31}, 
\beqn 
\min_{a_n \in \pi(x_n)} \rho_{n-1}(c_n(x_n,a_n))
\eeqn 
is l.s.c. in $x_n$ for any $x_n \in X_n$ $\bbP^{\pi}_{n-1}$-a.s. for $n \geq 1$.  Furthermore, by Lemma \ref{lem31}, we know that there exists an $\pi^* \in \Pi$ such that 
\beal 
\min_{a_n \in \pi(x_n)} \rho_{n-1}(c_n(x_n,a_n)) &= \rho_{n-1}(c_n(x_n,\pi^*(x_n)))\\
&= \rho_{n-1}(c_n(F_{n-1}(x_{n-1},a_{n-1},\xi_{n-1}),\\
&\indeq \pi^*(F_{n-1}(x_{n-1},a_{n-1},\xi_{n-1})))),\\
\eal
where $F_{n-1}$ is as defined in Equation \eqref{eqn23}, but we know that $\rho_{n-1}(c_n(x_n,\pi^*_n)$ is $\F_{n-1}$ measurable. Hence, the result follows by Lemma \ref{lem52}.
\eproof

For every $n \geq 0$, let $L_n(X_n)$ and $L_n(X_n,A_n)$ be the family of non-negative mappings on $(X_n,A_n)$, respectively. Denote
\begin{equation}
\label{eqn41}
T_{n}(v_{n+1}) := \min_{a_n \in A(x_n)} \big\{ c_n(x_n,a_n) + \gamma\rho_n(v_{n+1}(F_n(x_n,a_n, \xi_n)))\big\}.
\end{equation}
\belem \label{lem53}
Suppose that Assumption \ref{ass31}, \ref{ass32} and \ref{ass41} hold, then for every $n \geq 0$, we have
\beitem
\item[(a)] $T_n$ maps $L_{n+1}(X_{n+1})$ into $L_{n}(X_{n})$.
\item[(b)] For every $v_{n+1} \in L_{n+1}(X_{n+1})$, there exists a policy $\pi^*_n$ such that for any $x_n \in X_n$, $\pi^*_n(x_n) \in A_n(x_n)$ attains the minimum in \eqref{eqn41}, namely
\begin{equation}
\label{eqn316}
T_{n}(v_{n+1}) := c_n(x_n,\pi^*_n) + \gamma \rho(v_{n+1}(F_n(x_n,\pi^*_n, \xi_n)))
\end{equation}
\eitem
\elem
\beproof
By assumption, our one-step cost functions $c_n(x_n,a_n)$ are in $L_n(X_n)$. By Corollary \ref{cor51}, $\gamma \rho_n(v_{n+1}(F_n(x_n,\pi^*_n, \xi_n)))$ is in $L_n(X_n)$. Hence their sum is in $L_n(X_n,A_n)$, as well. Hence, the result follows via Corollary \ref{cor51} again.
\eproof
By Lemma \ref{lem53}, we express the optimality equations \eqref{eqn41} as 
\beqn
V_n^* = T_nV^*_{n+1}\; \trm{ for }n \geq 0.
\eeqn
Next, we continue with the following lemma. 
\belem 
\label{lem54}
Under the Assumptions \ref{ass31} and \ref{ass32}, for $n \geq 0$, let $v_n \in L_n(X_n)$ and $v_{n+1} \in L_{n+1}(X_{n+1})$.
\beitem
\item[(a)] If $v_n \geq T_n (v_{n+1})$, then $v_n \geq V_n^*$.
\item[(b)] If $v_n \leq T_n (v_{n+1})$ and in addition, 
\beqn
\lim_{N\rar \infty}v_{N}(x_{N+1}(\w)) = 0,
\eeqn
$\bbP$-a.s., then $v_n \leq V_n^*$.
\eitem
\elem 
\beproof
\beitem
\item[(a)]
By Lemma \ref{lem53}, there exists a policy $\pi=(\pi_n)_{n \geq 0}$ such that for all $n \geq 0$, 
\begin{equation}
v_{n}(x_n) \geq c_n(x_n,\pi_n) +  \rho_n( v_{n+1}(F_n(x_n,\pi_n, \xi_n))).
\end{equation}
By iterating the right hand side and by monotonicity of $\varrho_n(\cdot)$, we get 
\begin{equation}
v_{n}(x_{n}) \geq c_n(x_n,\pi_n) +  \varrho_n(\sum_{i = n+1}^{N-1}c_i(x_i,\pi_i) + v_{N}(x_{N})).
\end{equation}
Since $v_{N}(x_{N}) \geq 0$, we have 
\begin{equation}
v_{n}(x_{n}) \geq c_n(x_n,\pi_n) + \varrho_n(\sum_{i=n+1}^{N-1}c_i(x_i,\pi_i)), \trm{ a.s.}
\end{equation}
Hence, letting $N \rar \infty$, we obtain $v_{n}(x) \geq V_n(x,\pi)$ and so $v_{n}(x) \geq V_n^*(x)$.
\item[(b)] Suppose that $v_{n} \leq T_nv_{n+1}$ for $n \geq 0$, so that
\beqn
v_{n}(x_n) \leq  c_n(x_n,\pi_n) + \rho_n(c_{n+1}(x_{n+1},\pi_{n+1}) + v_{n+1}(x_{n+1}))
\eeqn
for any $\pi \in \Pi_{\trm{ad}}$, $\bbP^\pi_n$-a,s. Summing from $i=1$ to $i=N-1$ gives
\beal
&v_{n}(x_n) \leq c_n(x_n,a_n) +  \varrho_n(\sum_{i=1}^{N-1}c_{n+i}(x_{n+i},a_{n+i}) \\
&\indeq + \varrho_{N}(\sum_{i=n+N}^\infty c_i(x_i,a_i)))
\eal
Letting $N\rar \infty$ and by $\pi \in \Pi_{\mathrm{ad}}$, we get that 
\beqn
\lim_{N\rar \infty}\varo_n(v_{n+N}) = 0
\eeqn
so that we have 
\beqn 
v_n(x_n) \leq V_{n}(x_n,\pi),
\eeqn 
Taking infimum, we have
\beqn 
v_n(x_n) \leq V_{n}^*(x_n)
\eeqn 
Thus, we conclude the proof. 
\eitem
\eproof
To further proceed, we need the following technical lemma.
\belem 
\label{lem34}\cite{HL}
For every $N > n \geq 0$, let $X_n,A_n$ be complete, seperable metric spaces and $\bbK_n:= \{(x_n,a_n):x_n\in X_n,a_n\in A_n\}$ with $w_n$ and $w_{n,N}$ be functions on $\bbK_n$ that are non-negative, l.s.c. and inf-compact on $\bbK_n$. If $w_{n,N} \uar w_n$ as $N \rar \infty$, then 
\beqn
\lim_{N\rar \infty}\min_{a_n \in A_n}w_{n,N}(x_n,a_n) = \min_{a_n \in A_n} w_n(x_n,a_n),
\eeqn
for all $x_n \in X$.
\elem
The next result gives the validity of the convergence of value iteration. 
\begin{theorem}\label{thm51} Suppose that Assumptions \ref{ass31} and \ref{ass32} are satisfied. Then, for every $n \geq 0$ and $x_n \in X_n$,
\beqn
V_{n,N}^*(x_n) \uar V_n^*(x_n)\;\bbP_\trm{-a.s. }N\rar\infty
\eeqn 
and $V_n^*(x_n)$ l.s.c. $\bbP$-a.s.
\end{theorem}
\beproof
We obtain $V_{n,N}^*$ by the usual dynamic programming. Indeed, let $J_{N+1}(x_{N+1}) \equiv 0$ for all $x_{N+1} \in X_{N+1}$ a.s. and going backwards in time for $n=N,N-1,\ldots$, let
\beqn
\label{eqn328}
J_{n}(x_{n}) := \inf_{a_n \in A(x_n)} c_n(x_n,a_n) + \rho_n(J_{n+1}(F_n(x_n,a_n, \xi_n))).
\eeqn
Since $J_{N+1}(\cdot) \equiv 0$ is l.s.c., by backward induction, $J_N$ is  l.s.c. $\bbP$-a.s. and $\F_N$-measurable. Moreover, by Corollary \ref{lem31}, for every $t = N-1,...,n$, there exists $\pi_t^N$ such that $\pi_t^N(x_t) \in A_t(x_t)$ attains the minimum in Equation \eqref{eqn328}. Hence $\{ \pi_{N-1}^N,...,\pi_n^{N} \}$ is an optimal policy. We note that $c_n(x_n,a_n)$ as well as $\rho_n(J_{n+1}(F_n(x_n,a_n, \xi_n)))$ is  l.s.c., $\F_n$ measurable, inf-compact and non-negative. Hence their sum preserves those properties. Furthermore, $J_n$ is the optimal $(N-n)$ cost by construction. Hence, $J_n(x) = V_{n,N}^*(x)$ and since $J_n(x)$ is l.s.c. so is $V_{n,N}^*(x_n)$ with 
\beqn 
\label{eqn329}
V_{n,N}^*(x_n) := \inf_{a_n \in A(x_n)} \bigg(c_n(x_n,a_n) + \rho_n(V_{n+1,N}^*(x_{n+1})) \bigg). 
\eeqn
By the non-negativity assumption on $c_n(\cdot,\cdot)$ for all $n \geq 0$, the sequence $N \rar V_{n,N}^*$ is non-decreasing and $V_{n,N}^*(x_n) \leq V_n^*(x_n)$, for every $x_n \in X_n$ and $N > n$. Hence, denoting 
\beqn 
v_n(x_n) := \sup_{N>n}V_{n.N}^*(x_n) \trm{ for all }x_n \in X_n.
\eeqn 
and $v_n$ being supremum of l.s.c. functions is itself l.s.c. $\bbP$-a.s. and $\F$-measurable. Letting $N\rar \infty$ in \eqref{eqn329} by Lemma \ref{lem34}, we have that 
\beqn 
\label{eqn410} 
v_{n}(x_n) := \inf_{a_n \in A(x_n)} \bigg( c_n(x_n,a_n) + \rho_n(V_{n+1}(x_{n+1})) \bigg)
\eeqn
for all $n \in \bbN$ and $x_n \in X_n$. Hence, $v_n$ are solutions of the optimality equations, $v_n = T_nv_{n+1}$, and so by Lemma \ref{lem53}, $v_n(x_n) \geq V_n^*(x_n)$. This gives $v_n(x) = V_n^*(x)$. Hence, $V_{n,N}^* \uar V_n^*$ and $V_n^*$ is l.s.c. 
\eproof
Now, we are ready to prove our main theorem.
\beproof[Proof of Theorem \ref{thm21}] 
\beitem
\item[(a)] By Theorem \ref{thm51}, the sequence $(V_n^*)_{n \geq 0}$ is a solution to the optimality equations. By Lemma \ref{lem53}, it is the minimal such solution. 
\eitem
\item[(b)] By Theorem \ref{thm51}, the functions $V_n^*$ are l.s.c. $\bbP$-a.s. and $\F_n$-measurable. Therefore,
\beqn
\label{eqn32}
c_n(x_n,\pi^*_n) + \rho_n(V^*_{n+1}(x_{n+1}))
\eeqn
is non-negative, l.s.c. $\bbP$-a.s., $\F_n$-measurable and inf-compact on $\bbK_n$ for any $a_n \in A_n$, for every $n \geq 0$. Thus, the existence of optimal policy $\pi_n^*$ follows from Lemma \ref{thm51}. Iterating Equation \eqref{eqn32} gives 
\begin{align}
V_n^*(x_n) &= c_n(x_n, \pi_t^*) +  \varrho_n\bigg(\sum_{t=n+1}^{N-1}c_t(x_t, \pi_t^*) + V_{N}^*(x_{N})\bigg)\\
&\geq V_{n,N}(x_n,\pi_n^*).
\end{align} 
Letting $N \rar \infty$, we conclude that $V_n^*(x) \geq V_n(x,\pi^*)$. But by definition of $V_n^*(x)$, we have $V_n^*(x) \leq V_n(x,\pi^*)$. Hence, $V_n^*(x) = V_n(x,\pi^*)$, and we conclude the proof. 
\eproof

\subsection{An $\epsilon$-Optimal Approximation to Optimal Value}
We note that our iterative scheme via validity of convergence of value iterations in Theorem \ref{thm21} is computationally not effective for large horizon $N$ problem, since we have to calculate the dynamic programming equations for each time horizon $n \leq N$. To overcome this difficulty, we propose the following methodology, which requires only one time calculation of dynamic programming equations of the optimal control problem and is able to give an $\epsilon$-optimal approximation to the original problem. 

By Assumption \ref{ass41}, we have after some $N_0$
\beqn
\label{eqn470}
\varrho_{N_0}(\sum_{n=N_0+1}^\infty c_n(x_n, a_n)) < \epsilon \;\bbP \trm{-a.s}.
\eeqn 
But, then this means for the theoretical optimal policy $(\pi^*_n)_{n\geq0}$, justified in Theorem \ref{thm21}, we have 
\beqn
\varrho_{N_0}(\sum_{n=N_0+1}^\infty c_n(x_n, \pi^*_n)) \leq \epsilon\;\bbP \trm{-a.s.}
\eeqn
since, the optimal policy gives a smaller value than the one in Equation \eqref{eqn470}. 
Then, by monotonicity of $\varo$, for the optimal policy $\pi^*$ we have 
\beqn 
\varo_{N_0}(\sum_{n=N_0+1}^\infty c_n(x_n, \pi^*_n)) \leq \varo_{N_0}(\sum_{n=N_0+1}^\infty c_n(x_n, \pi_n)) \leq  \epsilon\;\bbP \trm{-a.s}.
\eeqn 
Hence, this means that by solving the optimal control problem up to time $N_0$ via dynamic programming and combine these decision rules $(\pi^*_0, \pi^*_1,\pi^*_2,...,\pi^*_{N_0})$ with the decision rules from time $N_0+1$ onwards, we have an $\epsilon$-optimal policy.  
Hence, we have proved the following theorem. 
\bethm\label{thm42}
Suppose that Assumptions \ref{ass31} and \ref{ass32} hold. Let $\pi_0 \in \Pi_{\trm{ad}}$ be the policy in Assumption \ref{ass41} such that 
\beqn
\varrho_{N_0}\big(\sum_{n=N_0+1}^\infty c_n(x_n, a_n)\big) < \epsilon\;\bbP \trm{-a.s}..
\eeqn 
Then, we have for the optimal policy
\beqn
\varrho_{N_0}\big(\sum_{n=N_0+1}^\infty c_n(x_n, a^*_n)\big) \leq \epsilon\;\bbP \trm{-a.s}.
\eeqn
Hence $\pi^* = \{\pi^*_0, \pi^*_1,\pi^*_2,...,\pi^*_{N_0}, \pi^0_{N+1},\pi^0_{N+2},\pi^0_{N+3}\dots \}$ is an $\epsilon$-optimal policy for the original problem. 
\ethm
\section{Applications}
\subsection{An Optimal Investment Problem}
In this section, we are going to study a variant of mean-variance utility optimization  (see e.g. \cite{BMZ}). The framework is as follows. We consider a financial market on an infinite time horizon $[0,\infty)$.  The market consists of a risky asset $S_n$ and a riskless asset $R_n$, whose dynamics are given by
\begin{align*}
&S_{n+1} - S_n = \mu S_n + \sigma S_n  \xi_n\\
&R_{n+1} - R_n = r R_n\\
\end{align*}
with $R_0 = 1, S_0=s_0$, where $(\xi_n)_{n\geq 0}$ are i.i.d standard normal random variables having distribution functions $\Phi$ on $\bbR$ with $\Z = L^1(\bbR,\cB(\bbR),\Phi)$ and $\mu,r,\sigma > 0$. We consider a self-financing portfolio composed of $S$ and $R$. We let $(\tilde{\pi}_n)_{n \geq 0}$ denote the amount of money invested in risky asset $S_n$ at time $n$ and $X_n$ denote the investor's wealth at time $n$. Namely, 
\beal 
X^{\tilde{\pi}}_n &= \tilde{\pi}_n S_n + R_n \\
X^{\tilde{\pi}}_{n+1}-X^{\tilde{\pi}}_n &= \tilde{\pi}_n(S_{n+1} - S_n) + (X^{\tilde{\pi}}_n - \tilde{\pi}_n) r R_n
\eal 
For each $n \geq 0$, we denote $\tilde{\pi}_n = X^{\tilde{\pi}}_n \pi_n$ so that $\pi_n$ stands for the fraction of wealth that is put in risky asset. Hence, the wealth dynamics are governed by 
\begin{align}
X^\pi_{n+1}-X^\pi_n &= [rZ^\pi_n + (\mu - r)\pi_n] + \sigma \pi_n \xi_n
\end{align}
with initial value $x_0 = S_0 + B_0$. We further assume $|\pi_n| \leq C$ for some constant $C > 0$ at each time $n \geq 0$.

The particular coherent risk measure used in this example is the mean-deviation risk measure that is in \textit{static setting} defined on $\Z$ as 
\begin{equation}
\varo(X) := \bbE^{\bbP}[X]  + \gamma g(X),
\end{equation} with $\gamma > 0$ with 
\beqn
g(X) :=  \bbE^{\bbP}\big(|X - \bbE^{\bbP}[X]| \big),
\eeqn 
for $X \in \Z$, where $\bbE^{\bbP}$ stands for the expectation taken with respect to the measure $\bbP$.
Hence $\gamma$ determines our \textit{risk averseness} level. For $\varo$ to satisfy the properties of a coherent risk measure, it is necessary that $\gamma$ is in $[0,1/2]$. In fact, $\gamma$ being in $[0,1/2]$ is both necessary and sufficient for $\varo$ to satisfy monotonicity (see  \cite{key-14}). Hence, for fixed $0 \leq \gamma \leq 1/2$ with $X \in \Z$, we have that 
\begin{equation}
\rho(X) = \sup_{m \in \cA}\la m, X\ra,
\end{equation}
where  $\cA$ is a subset of the probability measures, that are of the form (identifying them with their corresponding densities)
\beal
\label{eqn662}
&\cA = \bbpara m \in L^\infty(\bbR, \B(\bbR), \Phi): \int_\bbR m(x) d\Phi(x) = 1,\\
&\indeq m(x) = 1 + h(x) - \int_\bbR h(x)d\Phi(x),\; \| h \|_\infty \leq \gamma\;\Phi\trm{-a.s.} \bbpark
\eal
for some $h \in L^\infty(\bbR, \B(\bbR), \Phi)$.
Then, we \textit{define} for each time $n \geq 0$, the dynamic correspondent of $\rho$ as $\rho_n: \Z_{n+1} \rar \Z_n$ with
\begin{align}
\rho_n (X_{n+1}) &= \sup_{m_n \in \cA_{n+1} }\lan m_n, X_{n+1}\ran,
\end{align}
as in Equation $\eqref{eqn27},\eqref{eqn28},\eqref{eqn2140}$ using $(\bbR,\B(\bbR),\Phi)$. Hence, the controlled one step conditional risk mapping has the following representation 
\beqn
\sup_{m_n \in \cA_{n+1} } \lan m_n,X^\pi_n \ran,
\eeqn 
and our optimization problem reads as 
\beqn \label{eqn567}
\min_{\pi_n \in \Pi_{\mathrm{ad}}} \sup_{m_n \in \cA_{n+1} } \lan m_n,X^\pi_n \ran,
\eeqn 
where $\cA_{n+1}$ are the sets of conditional probabilities analogous to Equation \eqref{eqn662} with $\Pi_{\trm{ad}}$ as defined in Definition \ref{defn230}. Namely, $\cA_{n+1}$ is a subset of the conditional probability measures at time $n+1$ that are of the form (identifying them with their corresponding densities)
\beal
\label{eqn565}
&\cA_{n+1} = \bbpara m_{n+1} \in L^\infty(\Omega, \F_{n+1}, \bbP^\pi_{n+1}): \int_\O m_{n+1} d{\bbP}^\pi_{n+1} = 1,\\
&\indeq m_{n+1} = 1 + h - \int_\O h d{\bbP}^\pi_{n+1},\; \| h \|_\infty \leq \gamma\;\bbP^\pi_{n}\trm{-a.s.} \bbpark
\eal
for some $h \in L^\infty(\Omega, \F_{n+1}, \bbP^\pi_{n+1})$, where $\bbP^\pi_{n+1}$ stands for the conditional probability measure on $\O$ at time $n+1$ as constructed in \eqref{eqn2150}.  

Our one step cost functions are $c_n(x_n,a_n) = x_n$ for $n \geq 0$ for some discount factor $0 < \gamma < 1$ that are l.s.c. (in fact continuous) in $(x_n,a_n)$ for $n \geq 0$. 
Hence, starting with initial wealth at time 0, denoted by $x_0$, investor's control problem reads as 
\beal 
\label{exprob}
&x_0 + \min_{\pi \in \Pi_{\mathrm{ad}}} \varrho_0\bigg( \sum_{n = 1}^\infty  X^\pi_n \bigg) \\
&\triangleq x_0 + \min_{\pi \in \Pi_{\mathrm{ad}}} \lim_{N\rar \infty} \bigg( c_0(x_0,a_0) + \gamma\rho_0(c_1,(x_1,a_1) + \ldots+\gamma \rho_{N-1}(c_N(x_N,a_N) )\ldots)\bigg)
\eal 
We note that $\Pi_{\mathrm{ad}}$ is not empty so that our example satisfies Assumption \ref{ass41}. Indeed, by choosing $a_n \equiv 0$ for $n \geq 0$, i.e. investing all the current wealth into riskless asset $R_n$ for $n\geq 0$, we have that 
\beqn 
\varrho\bigg( \sum_{n = 0}^\infty  \gamma^n x_0 \bigg) = \frac{x_0}{1 - \gamma}
\eeqn 
Hence, as in Theorem \ref{thm42}, we find $N_0$ such that 
\beqn 
x_0 \sum_{n = N_0}^\infty \gamma^n < \epsilon.
\eeqn 
Thus, we write the corresponding \textit{robust} dynamic programming equations as follows. Starting with $V^*_{N_0+1} \equiv 0$ for $n = 1,2,...,N_0$, we have by Equation \eqref{eqn567}
\begin{align}
V^*_{n}(X^\pi_{n}) &= \min_{|\pi_n| \leq C} X^\pi_n + \gamma \rho_{n}(V^*_{n+1}(X^\pi_{n+1})) \\
&= \min_{|\pi_n| \leq C} X_n^\pi + \gamma \sup_{m_{n+1} \in \cA_{n+1}}\lan m_{n}, V^*_{n+1}(X_{n+1}^\pi)  \ran
\end{align}
going backwards iteratively at first stage, the problem to solve is then
\begin{align}
V^*_0(x_0) &= \min_{ |a_0| \leq C}x_0 +  \gamma\rho_0(  V^*_1( X^\pi_1) )\\
&= x_0 + \gamma \min_{|a_0| \leq C}\sup_{m_1 \in \cA_{1}}\lan m_1, V^*_1( X_1^\pi) \ran
\end{align}
Hence, the corresponding policy
\beqn 
\tilde{\pi} = \{ \pi^*_0, \pi^*_1,\pi^*_2,\ldots,\pi^*_{N_0},0,0,0,\ldots, \}
\eeqn
is $\epsilon$-optimal with the optimal value $V^\pi_0(x_0)$ for our example optimization problem \eqref{exprob}.

\subsection{The Discounted LQ-Problem}
We consider the linear-quadratic regulator problem in infinite horizon. We refer the reader to \cite{HL} for its study using expectation performance criteria. Instead of the expected value, we use the $\avr$ operator to evaluate total discounted performance.

For $n\geq 0$, we consider the scalar, linear system 
\beqn 
x_{n+1} = x_n + a_n + \xi_n,
\eeqn
with $X_0 = x_0$, where the disturbances $(\xi_n)_{n \geq 0}$ are independent, identically distributed random variables on $\Z_n^2 = L^2(\bbR,\B(\bbR),\bbP^n)$ with mean zero and $\bbE^{\bbP^n}[\xi^2_n] < \infty$. The control problem reads as
\beal 
\label{orn2}
&x_0 + \min_{\pi \in \Pi_{\mathrm{ad}}} \varrho_0\bigg( \sum_{n = 1}^\infty x^\pi_n \bigg) \\
&\triangleq x_0 + \min_{\pi \in \Pi_{\mathrm{ad}}} \lim_{N\rar \infty} \bigg( (x^2_0+ a^2_0) + \gamma\rho_0( (x^2_1+ a^2_1) \\
&\indeq + \ldots +\gamma^N\rho_{N-1}((x^2_N+ a^2_N) )\ldots)\bigg),
\eal 
where $\rho_n(\cdot): \Z_{n+1}^2 \rar \Z_{n}^2$ is the dynamic $\avr_\alf: \Z_{n+1}^2 \rar \Z_{n}^2$ operator defined as
\beal
\rho_n(Z) &\triangleq \sup_{m_{n+1} \in \cA_{n+1}} \lan m_{n+1},Z \ran,
\eal 
with
\beal
\label{eqn05999}
\cA_n = \big\{ &m_n \in L^{\infty}(\O,\F_n,\bbP^\pi_n): \int_\O m_n d\bbP^\pi_n = 1, \\
&\indeq  0 \leq \lVert m_n \rVert_\infty \leq \frac{1}{\alf}, \bbP^\pi_{n-1}-\trm{a.s.} \big\}
\eal 
We note that $\Pi_{\trm{ad}}$ is not empty. Indeed, choose $\pi_n \equiv 0$ for $n \geq 0$ so that 
\beqn 
x_n = x_0 + \sum_{i=0}^{n-1}\xi_i,
\eeqn 
with 
\beal
\varrho(\sum_{n=0}^\infty x^2_n  ) &\leq 2x_0^2 +  2\varrho( \sum_{n=0}^\infty \xi^2_n  ) \\
&\leq 2x^2_0 + 2\sum_{n=0}^\infty \gamma^n \avr_\alf(\xi^2_n)\\
&\leq 2 x^2_0 + 2 \sum_{n=0}^\infty \gamma^n \frac{1}{\alf}\bbE^{\bbP}[\xi^2_i]\\
&\leq 2 x^2_0 + \frac{2\sg^2}{\alf(1-\gamma)}\\
&< \infty,
\eal 
where we used Equation \eqref{eqn05999} in the third inequality. 
Hence, we find $N_0$ such that 
\beqn 
\frac{2\sg^2}{\alf}\sum_{n=N_0}^\infty \gamma^n < \epsilon.
\eeqn 
Starting with $J_{N_0+1} \equiv 0$, the corresponding $\epsilon$-optimal policy for $n=0,1,\ldots,N_0$ is found via 
\beal
J_n(x_n) = \min_{|\pi_n|\leq C}\bigg( (x^2_n + a^2_n) + \gamma \avr^n_{\alf}\big(J_{n+1}(x_{n+1} + a_{n+1})\big) \bigg),
\eal 
so that at the final stage, we have
\beal
J_0(x_0) &= \min_{ |a_0| \leq C}x_0 +  \gamma\avr_\alf(  J_1( x^\pi_1) )\\
&= x_0 + \gamma \min_{|a_0| \leq C}\sup_{m_1 \in \cA_{1}}\lan m_1, J_1( x_1) \ran,
\eal
where $\cA$ is as defined in Equation \eqref{eqn05999}. Thus, the corresponding policy
\beqn 
\tilde{\pi} = \{ \pi^*_0, \pi^*_1,\pi^*_2,\ldots,\pi^*_{N_0},0,0,0,\ldots, \}
\eeqn
is $\epsilon$-optimal with the optimal value $V^\pi_0(x_0)$ for problem \eqref{orn2}.

\end{document}